\documentclass[twoside, 11pt]{article}

\usepackage{mathrsfs,amsfonts,amsmath}
\RequirePackage{ifthen,calc}
\usepackage{theorem}
\usepackage[dvips]{color}

\def\R{{\mathbb R}}

\def\E{{\mathbb E}}
\def\P{{\mathbb P}}

\def\N{{\mathbb N}}

 \catcode`@=11
 \def\@evenhead{\hbox to\textwidth{\footnotesize\rm\thepage \hfill
  {\it }}} 

 \def\@oddhead{\hbox to \textwidth{\footnotesize{\it
  OFBMS and Martingale Differences } \hfill\thepage}}

 \renewcommand{\section}{\makeatletter
 \renewcommand{\@seccntformat}[1]{{\csname the##1\endcsname.}\hspace{0.45em}}
 \makeatother \@startsection
{section}
{1}
{0pt}
{\baselineskip}
{0.5\baselineskip}
{\normalsize\bfseries\mathversion{bold}}}

\renewcommand{\subsection}{\makeatletter
 \renewcommand{\@seccntformat}[1]{{\csname the##1\endcsname.}\hspace{0.45em}}
 \makeatother \@startsection
{subsection}
{1}
{0pt}
{\baselineskip}
{0.5\baselineskip}
{\normalsize\bfseries\mathversion{bold}}}

\catcode`@=12

\newtheorem{theorem}{\noindent Theorem}[section]
\newtheorem{lem}{\noindent Lemma}[section]

\theoremheaderfont{\normalfont\bfseries}
\theorembodyfont{\slshape} {\theorembodyfont{\rmfamily}} {\theorembodyfont{\rmfamily}
\newtheorem{defn}{\noindent Definition}[section]}
{\theorembodyfont{\rmfamily} } {\theorembodyfont{\rmfamily}
}
{\theorembodyfont{\rmfamily} } {\theorembodyfont{\rmfamily}
\newtheorem{rem}{\noindent Remark}[section]}
{\theorembodyfont{\rmfamily} } {\theorembodyfont{\rmfamily} } {\theorembodyfont{\rmfamily}
}

 \def\beqlb{\begin{eqnarray}}\def\eeqlb{\end{eqnarray}}
 \def\beqnn{\begin{eqnarray*}}\def\eeqnn{\end{eqnarray*}}
 \numberwithin{equation}{section}
\def\2R{\mathbb{R}_+\times\mathbb{R}}
\def\qed{\hfill$\square$\smallskip}

\def\3R{\mathbb{R}_+\times\mathbb{R}_-}

\begin{document}
\title{\bf OPERATOR FRACTIONAL  BROWNIAN SHEET AND MARTINGALE DIFFERENCES
\footnotetext{\hspace{-5ex}
${[1]}$ School of Statistics, Shandong University of Finance and Economics,
Jinan, 250014  China
\\${[2]}$  Department of Mathematics, Anhui Normal  University,
 Wuhu,  241000  China
 \\ ${\S}$ Corresponding author.
 \newline
}}
\author{\small Hongshuai Dai$^{1}$, Guangjun Shen ${}^{\S}$ $^2$ and Liangwen Xia $^2$ }

\maketitle

\begin{abstract}
In this paper, inspired by the fractional Brownian sheet of Riemann-Liouville type, we introduce the operator fractional Brownian sheet of Riemman-Liouville type, and  study some properties of it. We also present an approximation in law to it  based on the martingale differences.
\end{abstract}


\vspace{2mm}

\section{Introduction}
Self-similar processes, first studied rigorously by Lamperti \cite{Lamperti} under the name ``semi-stable",
are stochastic processes that are invariant in distribution under suitable scaling of time and space.
There has been an extensive literature on self-similar processes. We refer to Vervaat \cite{Vervaat} for
general properties, to Samorodnitsky and Taqqu \cite{Samorodnitsky and Taqqu}[Chaps.7 and 8] for studies
on Gaussian and stable self-similar processes and random fields.

The fractional Brownian motion (fBm) as a well-known self-similar process has been studied extensively.
Many results about weak approximation to fBms  have been established recently.
See \cite{Yu 2, Yu 1} and the references therein.
 We point out that the fBm does not represent a casual time-invariant system as there
is no well-defined impulse response function.  Hence, based on the  Riemann-Liouville fractional integral, Barnes and Allan \cite{BA1966} introduced the
fractional Riemann-Liouville (RL) Brownian motion (RL-fBm). RL-fBms share with fBms many properties which include self-similarity, regularity of
sample paths, etc.- with one notable exception that its increment process is nonstationary.  For more information on RL-fBms, refer to Lim \cite{Lim2001} and the references therein.  On the other hand, there are two typical multiparameter extensions of fBms, one of which is
the fractional Brownian sheet introduced by Kamont \cite{Kamont}. Fractional Brownian
sheets have been studied extensively as a representative of anisotropic Gaussian random
fields. For more information,
refer to   \cite{XAVIER 1, XAVIER 2} and \cite{Zhi Wang 1, Zhi Wang 2}.  Inspired by the study of RL-fBms and fractional Brownian sheets,  Dai \cite{Dai2015} introduced the multifractional Riemann-Liouville Brownian sheet and studied the weak limit theorem for it.

The definition of self-similarity has been extended to allow scaling by linear operators on multidimensional space $\mathbb{R}^d$,
and the corresponding processes are called operator self-similar processes.
We refer to \cite{Laha}, \cite{Lamperti}, \cite{Mason and Xiao} and the references therein.
We note that Didier and Pipiras  \cite{Didier and Pipiras 1, Didier and Pipiras 2} introduced
the operator fractional Brownian motions (ofBm in short) as an extension of fBms and studied their properties.
Similar to fBms, weak limit theorems for ofBms have also attracted a lot of interest.
Recently, Dai and his  coauthors \cite{Hong Shuai Dai 2}-\cite{Hong Shuai Dai 1} presented  some weak limit theorems for some kinds of  ofBms.

In contrast to the extensive study on the multiparameter extension of fBms,
there is little work studying the multiparameter extension of ofBms. Inspired by the study of the fractional Brownian sheet and the operator fractional Brownian motion of Riemann-Liouville type introduced by Dai \cite{Hong Shuai Dai 3}, we will introduce a new  random field,
which we call  the operator fractional Brownian sheet of Riemann-Liouville type, and present an approximation to it.

Most of the estimates of this paper contain unspecified constants.
An unspecified positive and finite constant will be denoted by $C$,
which may not be the same in each occurrence. Sometimes we shall
emphasize the dependence of these constants upon parameters.

At the end of this section, we point out that all processes  considered here
are assumed to be proper.  We say that a process $\{X(t); t \in \R^d_+\} $ is proper if for each $t\in\R_+^d$ the distribution of $X(t)$ is full; that is, the
distribution is not contained in a proper hyperplane.

The rest of this paper is organized as follows. In Section 2, we  introduce the operator fractional Brownian sheet of Riemann-Liouville type and state some properties. We present an approximation in law to it in Section 3.  A final note is presented at the end of this paper.

\section{Operator Fractional Brownian Sheet}\label{sec2}
In this section, we first introduce the operator fractional Brownian sheet of Riemann-Liouville type and then study some properties of it. For any $x\in\mathbb{R}^d$, $x^{T}$ denotes the transpose of $x$.
Let $(\Omega, \mathcal{F}, \mathbb{P})$ be a probability space and  $\{\mathcal{F}_{t,s};\, (t,s)^T\in\R_+^2\}$ be a family of sub-$\sigma$-fields of $\mathcal{F}$ such that $\mathcal{F}_{t,s}\subseteq\mathcal{F}_{t', s'}$ for any $(t, s)^T<(t', s')^T$ with the usual partial order.   Moreover, for any stochastic process $Y=\{Y(t,s); \, (t,s)^T\in\R_+^2\}$, we denote by $\Delta_{(t, s)}Y(t^{'}, s^{'})$ the increment of $Y$ over the
rectangle $(t, \,t']\times (s,\, s^{'}]$, that is,
$$\Delta_{(t,\; s)}Y(t^{'}, s^{'})=Y(t^{'}, s^{'})-Y(t, s^{'})-Y(t^{'}, s)+Y(t, s).$$
Let $\sigma(A)$ be the collection of all eigenvalues of a linear operator $A$ on $\R^d$. Let
\beqnn
\lambda_{A}=\min\{\textrm{Re}\lambda: \lambda\in\sigma(A)\} ~~\textrm{and} ~~\Lambda_{A}=\max\{\textrm{Re}\lambda: \lambda\in\sigma(A)\}.
\eeqnn
Moreover,  given any linear operator $A$ on $\R^d$ and $t>0$, we define the power operator
\beqnn
t^A=\sum_{k=0}^\infty (\log t)^k \frac{A^k}{k!}.
\eeqnn
Next,
we recall the operator fractional Brownian motion of Riemann-Liouville type introduced by  Dai \cite{Hong Shuai Dai 3}.
Let $D$ be a linear operator on $\mathbb{R}^{d}$ with $0<\lambda_{D}, \Lambda_{D}<1$.
We define the operator fractional Brownian motion of Riemann-Liouville type $\tilde{X}=\{\tilde{X}(t); t\in\R_+\}$ with exponent $D$ by
\begin{align}\label{da-36}
\tilde{X}(t)=\int_{0}^{t}(t-u)^{D-I/2}dW(u),
\end{align}
where $W(u)=\{W^{1}(u), ... , W^{d}(u)\}^{T}$ is a standard $d$-dimensional Brownian motion and $I$ is the $d\times d$ identity matrix.

Based on \eqref{da-36}, we can define the operator fractional Brownian sheet of Riemann-Liouville type $X=\big\{X(t,s); (t,s)^T\in \R_+^2\big\}$ as follows.

\begin{defn}\label{da-35} Let  $\tilde{B}$ be the standard Brownian sheet.
The operator fractional Brownian  sheet of Riemann-Liouville type $X=\big\{X(t,s); (t,s)^T\in\R_+^2\big\}$ is defined by
\begin{align}\label{da-8}
X(t, s)=\int_{0}^{t}\int_{0}^{s}(t-u)^{\frac{D}{2}-\frac{I}{2}}(s-v)^{\frac{D}{2}-\frac{I}{2}}B(du, dv),
\end{align}
where $B(du, dv)=\big(B^{1}(du, dv), ..., B^{d}(du, dv)\big)^T$ with $B^i$ being independent copies of $\tilde{B}$, and $D$ is a linear operator on $\R^d$ with $0<\lambda_D, \Lambda_D<1$.
\end{defn}

\begin{rem} Let $x_+ = \max\{x, 0\}$.
From \eqref{da-8} and Mason and Xiao \cite{MM1994}, we get that  $X$  is an $\R^d$-valued Gaussian random field with mean zero vector and  for any $(t_1,t_2)^T, (s_1,s_2)^T\in\R_+^2$
\beqlb\label{da-33}
\E\big[X(t_1,t_2)X^T(s_1,s_2)\big]=&&\int_0^\infty\int_0^\infty \big[(t_1-u)_+(t_2-v)_+]^{\frac{D}{2}-\frac{I}{2}}\nonumber
\\&&\quad\cdot\big[(s_2-v)_+(s_1-u)_+\big]^{\frac{D^*}{2}-\frac{I}{2}}dudv,
\eeqlb
where $D^*$ is the adjoint operator of $D$.
\end{rem}

It is obvious that the equation \eqref{da-8} is well defined. Next, we study some properties of the random field $X$. We first introduce the following notation.  Let $\|x\|_2$ denote the usual Euclidean norm of $x\in \R^d$. Similar to  Dai, Shen and Kong \cite{DSK},
 $End (\R^d)$  denotes  the set of linear operators on $\R^d$
(endomorphisms). Furthermore, we will not distinguish an operator $D\in End (\R^d)$ from its associated
matrix relative to the standard basis of $\R^d$.  For any $A\in End (\R^d)$, let $\|A\| = \max_{\|x\|_2=1}\|Ax\|_2$ be the operator norm of $A$. Next, we  recall the definition of operator self-similar processes.   Recall that an $\R^d$-valued stochastic  process $\tilde{Y}=\{\tilde{Y}(t); t \in \R^2_+\}$ is said to be operator self-similar (o.s.s.) if
it is  continuous in law at each $t\in \R_+^2$, and there exists $D\in End(\R^d)$  such that
\beqnn
\big\{\tilde{Y}(ct)\big \}\stackrel{\mathscr{D}}{=}\big\{c^D \tilde{Y}(t)\big\}\;\textrm{for all }\;c>0,
\eeqnn
where $\stackrel{\mathscr{D}}{=}$ denotes the equality of all finite-dimensional distributions.

\begin{theorem}
The random field $X=\{X(t,s);\;(t,\,s)^T\in\R_+^2\}$ is an operator self-similar Gaussian random field with exponent $D$. Moreover, $X$
has a version with continuous sample paths a.s..
\end{theorem}
{\it Proof:}  We first check the operator self-similarity.  For every $c>0$, we have
\beqnn
X(ct,cs)&&=\int_0^\infty\int_0^\infty(ct-u)_+^{\frac{D}{2}-\frac{I}{2}}(cs-v)_+^{\frac{D}{2}-\frac{I}{2}}dB(u,v)
\\&&\stackrel{\mathscr{D}}{=}c^{D-I}\int_0^{ct}\int_0^{cs}(t-\frac{u}{c})^{\frac{D}{2}-\frac{I}{2}}(s-\frac{v}{c})^{\frac{D}{2}-\frac{I}{2}}dB(u,v)
\\&&\stackrel{\mathscr{D}}{=}c^D X(t,s),
\eeqnn
since \beqlb\label{da-26}B(cu,\;cv)\stackrel{\mathscr{D}}{=}c^{I}B(u,\;v)\;\textrm{and}\;z^Dy^D=(zy)^D\;\textrm{for any}\;z>0,y>0.\eeqlb

Next, we check the sample continuity.   Choose any $t=(t_1,t_2)^T,s=(s_1,s_2)^T\in\R_+^2$.  Without loss of generality, we assume that $s<t$ with the usual partial order, and  $\|t-s\|_2\leq 1$. By some calculations, we have
\beqlb
&&\Delta_{s}X(t)=\nonumber
\\&&\quad\int_0^\infty\int_0^\infty\Big(\big(t_1-u\big)_+^{\frac{D}{2}-\frac{I}{2}}-\big(s_1-u\big)_+^{\frac{D}{2}-\frac{I}{2}}\Big)
\Big(\big(t_2-v\big)_+^{\frac{D}{2}-\frac{I}{2}}-\big(s_2-v\big)_+^{\frac{D}{2}-\frac{I}{2}}\Big)B(du,dv).\qquad
\eeqlb
Hence,
\beqlb\label{da-1}
\big\|\Delta_sX(t)\big\|_2^2=\sum_{i=1}^d\Big(\int_0^\infty\int_0^\infty \sum_{j=1}^d F_{i,j}(t,s,u,v)B^j(du,dv)\Big)^2,
\eeqlb
where \beqnn
F(t,s,u,v)&&=\Big((t_1-u)_+^{\frac{D}{2}-\frac{I}{2}}-(s_1-u)_+^{\frac{D}{2}-\frac{I}{2}}\Big)
\Big((t_2-v)_+^{\frac{D}{2}-\frac{I}{2}}-(s_2-v)_+^{\frac{D}{2}-\frac{I}{2}}\Big)
\\&&=\big(F_{i,j}(t,s,u,v)\big)_{d\times d}.
\eeqnn
Noting that $\int_0^\infty\int_0^\infty \sum_{j=1}^d F_{i,j}(t,s,u,v)B^j(du,dv)$ is a Gaussian random variable, we get from \eqref{da-1} that for any even $k\in\N$
\beqlb\label{da-5}
\E\Big[\big\|\Delta_sX(t)\big\|_2\Big]^k\leq C\Big[\int_0^\infty\int_0^\infty\big\|F(t,s,u,v)\big\|^2dudv\Big]^{\frac{k}{2}}.
\eeqlb
On the other hand, we have
\beqlb\label{da-2}
\big\|F(t,s,u,v)\big\|^2\leq C\|F_1(t,s,u,v)\|^2 \times\|F_2(t,s,u,v)\|^2,
\eeqlb
where
\beqnn
F_1(t,s,u,v)=\big(t_1-u\big)_+^{\frac{D}{2}-\frac{I}{2}}-\big(s_1-u\big)_+^{\frac{D}{2}-\frac{I}{2}},
\eeqnn
and
\beqnn
F_2(t,s,u,v)=\big(t_2-v\big)_+^{\frac{D}{2}-\frac{I}{2}}-\big(s_2-v\big)_+^{\frac{D}{2}-\frac{I}{2}}.
\eeqnn
Now, we look at
\beqnn
\int_0^\infty\big\|F_1(t,s,u,v)\big\|^2dudv.
\eeqnn
 By using the same method as in Dai, Hu and Lee \cite{Hong Shuai Dai 1}, we have
\beqlb\label{da-6}
\int_0^\infty\|(t_1-u)_+^{\frac{D}{2}-\frac{I}{2}}-(s_1-u)_+^{\frac{D}{2}-\frac{I}{2}}\|^2du\leq C(t_1-s_1)^{\lambda_D-\delta}.
\eeqlb
Similarly,
\beqlb\label{da-28}
\int_0^\infty \big\|F_2(t,s,u,v)\big\|^2dv\leq C (t_2-s_2)^{\lambda_D-\delta}.
\eeqlb
From Maejima and Mason \cite{MM1994}, and \eqref{da-5}-\eqref{da-28}, we can get that for any $\delta>0$ with $\lambda_D-\delta>0$,
\beqlb\label{da-29}
\E\Big[\big\|\Delta_sX(t)\big\|_2^k\Big]&&\leq C\Big[(t_1-s_1)^{\lambda_D-\delta}\times(t_2-s_2)^{\lambda_{D}-\delta}\Big]^{\frac{k}{2}}\nonumber
\\&&\leq C\big\|t-s\big\|_2^{(\lambda_D-\delta)k}.
\eeqlb
The sample continuity  follows from Garsia \cite{Garsia} and \eqref{da-29}. \qed
\section{Limit Theorem }
 One  aim of this paper is to present an approximation in law to the operator fractional Brownian sheet of Riemann-Liouville type $X$ via the martingale differences.  In order to reach it, we first recall some facts about the martingale differences.  Similar to Wang, Yan and Yu \cite{Zhi Wang 1}, we  use the definitions and notations introduced in the basic work of Cairoli and Walsh \cite{Cairoli and Walsh}  on stochastic calculus in the plane.  For any ${\bf n}=(n_1,n_2)^T\in N_0\times M_0$ with $N_0=\{1,\cdots,n_0\}$ and $M_0=\{1,\cdots,m_0\}$, let
$\tilde{\mathcal{F}}_{\bf n}:=\mathcal{F}_{n_1,m_0}\bigvee\mathcal{F}_{n_0,n_2}$, the $\sigma-$ fields generated by $\mathcal{F}_{n_1,m_0}$ and
 $\mathcal{F}_{n_0,n_2}$.
Now, we recall the definition of the strong martingale.
\begin{defn}
An integrable process $Y=\{Y({\bf n}), {\bf n}\in N_0\times M_0\}$ is called a strong martingale if:
\begin{itemize}
\item[(i)] $Y$ is adapted;
\item[(ii)] $Y$ vanishes on the axes;
\item[(iii)]$\E\big[\Delta_{\bf n}Y({\bf m})|\tilde{\mathcal{F}}_{\bf n}\big]=0$ for any ${\bf n}\leq {\bf m}\in N_0\times M_0$ with the usual partial order.
\end{itemize}
\end{defn}
Let  $\big\{\xi^{(n)}=(\xi_{i,j}^{(n)}, \mathcal{F}_{i,j}^{(n)})\big\}_{n\in\N}$ be a sequence such that for all
$$
\E\Big[\xi_{i+1,j+1}^{(n)}|\mathcal{F}_{i,j}^{(n)}\Big]=0,
$$
where $\mathcal{F}_{i,j}^{(n)}=\mathcal{F}^{(n)}_{i,n}\bigvee\mathcal{F}^{(n)}_{n,j}$ with $\mathcal{F}_{k,l}^{(n)}$ being the $\sigma-$ fields generated by all $\xi_{r,s}^{(n)}$, $r\leq k,\;s\leq l$.
Then we call   $\big\{\xi^{(n)}=(\xi_{i,j}^{(n)}, \mathcal{F}_{i,j}^{(n)})\big\}_{n\in\N}$  a martingale differences sequence.

It is well known that if the martingale differences sequence $\{\xi^{(n)}\}$ satisfies the following condition
\beqnn
\sum_{i=1}^{\lfloor nt \rfloor}\sum_{j=1}^{\lfloor ns \rfloor} \big(\xi_{i,j}^{(n)}\big)^2\to t\cdot s
\eeqnn
in the sense of $\mathcal{L}^1$, then  the sequence
\beqnn
\sum_{i=1}^{\lfloor nt \rfloor}\sum_{j=1}^{\lfloor ns \rfloor} \xi_{i,j}^{(n)}
\eeqnn
converges weakly to the Brownian sheet, as $n$ goes to infinity (see for example, Morkvenas \cite{M1984}.)  Recently,  Wang, Yan and Yu \cite{Zhi Wang 1} extended this work to the fractional Brownian sheet. If  $\{\xi^{(n)}\}$ is a square integrable martingale differences sequence satisfying the following two conditions:
\begin{align}\label{a-1}
\lim_{n\rightarrow\infty}n(\xi_{i, j}^{(n)})=1, ~~a.s.
\end{align}
for any $1\leq i,j\leq n$,
and
\begin{align}\label{a-2}
\max_{1\leq i, j\leq n}|\xi_{i, j}^{(n)}|\leq\frac{C}{n}, a.s.
\end{align}
for some $C\geq1$, then, based on $\{\xi^{(n)}\}$, the authors of \cite{Zhi Wang 1} constructed a sequence  to converge weakly to the fractional Brownian sheet.   Inspired by these results, we want to study the weak limit theorem for the operator fractional Brownian sheet of Riemann-Liouville type $X$ introduced in Definition \ref{da-35}.  Similar to  Wang, Yan and Yu \cite{Zhi Wang 1}, we  assume that $\frac{1}{2}<\lambda_D,\Lambda_D<1$ in the rest of this paper.

Define
\begin{align}
\eta_{i, j}^{(n)}=(\xi_{i, j, 1}^{(n)}, ... , \xi_{i, j, d}^{(n)})^{T},
\end{align}
and
\begin{align}
B_{n}(t, s)=\sum_{i=1}^{\lfloor nt\rfloor}\sum_{j=1}^{\lfloor ns\rfloor}\eta_{i, j}^{(n)},
\end{align}
where $\xi_{i,j,k}^{(n)},k=1,\cdots,d,$ are independent copies of $\xi_{i,j}^{(n)}$.

From the above arguments, we obtain that $\big\{(\eta_{i, j}^{(n)}, \mathcal{F}_{i, j}^{(n)})\big\}_{n\in\mathbb{N}}$ is still a sequence of square integrable martingale
differences on the probability space $(\Omega, \mathcal{F}, \mathbb{P})$.

For any $n\geq1$ and $(t, s)^T\in[0,\,1]^2$, define
\begin{align}\label{da-30}
X_{n}(t, s)&=\int_{0}^{t}\int_{0}^{s}(t-u)_{+}^{\frac{D}{2}-\frac{I}{2}}(s-v)_{+}^{\frac{D}{2}-\frac{I}{2}}B_{n}(du, dv)\nonumber\\
&=n^{2}\sum_{i=1}^{\lfloor nt\rfloor}\sum_{j=1}^{\lfloor ns\rfloor}\eta_{i, j}^{(n)}
\int_{\frac{i-1}{n}}^{\frac{i}{n}}\int_{\frac{j-1}{n}}^{\frac{j}{n}}(t-u)_{+}^{\frac{D}{2}-\frac{I}{2}}(s-v)_{+}^{\frac{D}{2}-\frac{I}{2}}dudv.
\end{align}
Then, we have the following approximation. As a prelude to giving the result, let
\beqnn
\mathcal{D}([0,\;1]^2)&&=\mathcal{D}([0,\;1]^2, \R^d).
\eeqnn

\begin{theorem}\label{thm-2} Let $\frac{1}{2}<\lambda_D,\Lambda_D<1$.
The sequence of processes $\{X_{n}(t, s);(t, s)^T\in[0, 1]^2\}$ given by \eqref{da-30} converges weakly, as $n\to\infty$  in $\mathcal{D}([0,\;1]^2)$,
 to the operator fractional Brownian sheet of Riemann-Liouville type $\{X(t,s);(t,s)^T\in[0,\;1]^2\}$ given by \eqref{da-8}.
\end{theorem}

The proof of Theorem \ref{thm-2} is based on a series of technical results.
\begin{lem}\label{lem-2}
Let $\{X_{n}(t, s)\}$ be the family of processes defined by \eqref{da-30}. Then for any $s=(s_1,s_2)^T<t=(t_1, t_2)^T<u=(u_1, u_2)^T\in[0,\;1]^2$,
\begin{align}\label{da-34}
\E\Big[\big\|\Delta_{s}X_{n}(t)\big\|_2\big\|\Delta_{t}X_{n}(u)\big\|_2\Big]^2 \leq C(u_2-s_2)^{2H}(u_1-s_1)^{2H},
\end{align}
where $H=\lambda_D-\delta$ with $0<\delta<\lambda_D-\frac{1}{2}$.
\end{lem}
{\it Proof:} From \eqref{da-30}, we have
\begin{align*}
&\Delta_{s}X_{n}(t)
\\&=\int_{s_1}^{t_1}\int_{s_2}^{t_2}\Big((t_1-u)_{+}^{\frac{D}{2}-\frac{I}{2}}-(s_1-u)_{+}^{\frac{D}{2}-\frac{I}{2}}\Big)
\Big((t_2-v)_{+}^{\frac{D}{2}-\frac{I}{2}}-(s_2-v)_{+}^{\frac{D}{2}-\frac{I}{2}}\Big)B_{n}(du, dv)\nonumber\\
&=\sum_{i=1}^{\lfloor nt_1\rfloor}\sum_{j=1}^{\lfloor nt_2\rfloor}n^{2}\int_{\frac{i-1}{n}}^{\frac{i}{n}}\int_{\frac{j-1}{n}}^{\frac{j}{n}}
\Big((\frac{\lfloor nt_1\rfloor}{n}-u)_{+}^{\frac{D}{2}-\frac{I}{2}}-(\frac{\lfloor ns_1\rfloor}{n}-u)_{+}^{\frac{D}{2}-\frac{I}{2}}\Big)\nonumber\\
&\quad\times\Big((\frac{\lfloor nt_2\rfloor}{n}-v)_{+}^{\frac{D}{2}-\frac{I}{2}}-(\frac{\lfloor ns_2\rfloor}{n}-v)_{+}^{\frac{D}{2}-\frac{I}{2}}\Big)dudv\eta_{i, j}^{(n)}.
\end{align*}
It follows from \eqref{a-2} and the Cauchy-Schwarz inequality that
\beqlb\label{da-3}
\E\Big[\big\|\Delta_{s}X_{n}(t)\big\|_2\Big]^{4}&=\E\bigg[\Big\|\sum_{i=1}^{\lfloor nt_1\rfloor}\sum_{j=1}^{\lfloor nt_2\rfloor}n^{2}
\int_{\frac{i-1}{n}}^{\frac{i}{n}}\int_{\frac{j-1}{n}}^{\frac{j}{n}}\Big((\frac{\lfloor nt_1\rfloor}{n}-u)_{+}^{\frac{D}{2}-\frac{I}{2}}-(\frac{\lfloor ns_1\rfloor}{n}-u)_{+}^{\frac{D}{2}-\frac{I}{2}}\Big)\nonumber\\
&\quad\times\Big((\frac{\lfloor nt_2\rfloor}{n}-v)_{+}^{\frac{D}{2}-\frac{I}{2}}-(\frac{\lfloor ns_2\rfloor}{n}-v)_{+}^{\frac{D}{2}-\frac{I}{2}}\Big)dudv\eta_{i, j}^{(n)}\Big\|_2\bigg]^{4}\nonumber
\\
&\leq C\sum_{i=1}^{\lfloor nt_1\rfloor}\Big(\int_{\frac{i-1}{n}}^{\frac{i}{n}}\|(\frac{\lfloor nt_1\rfloor}{n}-u)_{+}^{\frac{D}{2}-\frac{I}{2}}-(\frac{\lfloor ns_1\rfloor}{n}-u)_{+}^{\frac{D}{2}-\frac{I}{2}}\|^{2}du\Big)^2\nonumber\\
&\quad\times \sum_{j=1}^{\lfloor nt_2\rfloor}\Big(\int_{\frac{j-1}{n}}^{\frac{j}{n}}\|(\frac{\lfloor nt_2\rfloor}{n}-v)_{+}^{\frac{D}{2}-\frac{I}{2}}-(\frac{\lfloor ns_2\rfloor}{n}-v)_{+}^{\frac{D}{2}-\frac{I}{2}}\|^{2}dv\Big)^2\nonumber\\
&\leq C\bigg(\int_{0}^{t_1}\|(\frac{\lfloor nt_1\rfloor}{n}-u)_{+}^{\frac{D}{2}-\frac{I}{2}}-(\frac{\lfloor ns_1\rfloor}{n}-u)_{+}^{\frac{D}{2}-\frac{I}{2}}\|^{2}du\bigg)^2\nonumber
\\&\quad\bigg(\int_{0}^{t_2}\|(\frac{\lfloor nt_2\rfloor}{n}-v)_{+}^{\frac{D}{2}-\frac{I}{2}}-(\frac{\lfloor ns_2\rfloor}{n}-v)_{+}^{\frac{D}{2}-\frac{I}{2}}\|^{2}dv\bigg)^2\nonumber\\
&\leq C\bigg(\int_{0}^{1}\|(\frac{\lfloor nt_1\rfloor}{n}-u)_{+}^{\frac{D}{2}-\frac{I}{2}}-(\frac{\lfloor ns_1\rfloor}{n}-u)_{+}^{\frac{D}{2}-\frac{I}{2}}\|^{2}du\bigg)^2\nonumber
\\ &\bigg(\int_{0}^{1}\|(\frac{\lfloor nt_2\rfloor}{n}-v)_{+}^{\frac{D}{2}-\frac{I}{2}}-(\frac{\lfloor ns_2\rfloor}{n}-v)_{+}^{\frac{D}{2}-\frac{I}{2}}\|^{2}dv\bigg)^2.
\eeqlb
From Dai, Hu and Lee \cite{Hong Shuai Dai 1}, we obtain that
\begin{align*}
\int_{0}^{1}\parallel(\frac{\lfloor nt_1\rfloor}{n}-u)_{+}^{\frac{D}{2}-\frac{I}{2}}-(\frac{\lfloor ns_1\rfloor}{n}-u)_{+}^{\frac{D}{2}-\frac{I}{2}}\parallel^{2}du\leq C(\frac{\lfloor nt_1\rfloor}{n}-\frac{\lfloor ns_1\rfloor}{n})^{H},
\end{align*}
where  $H=\lambda_{D}-\delta$. Then \eqref{da-3} can be bounded by
\begin{align}\label{da-4}
C\Big(\frac{\lfloor  nt_1 \rfloor-\lfloor ns_1 \rfloor}{n}\Big)^{2H}
\Big(\frac{\lfloor nt_2 \rfloor-\lfloor ns_2\rfloor}{n}\Big)^{2H}.
\end{align}
Hence, for  any $s<t<u\in[0,\;1]^2$, we have
\beqlb\label{da-27}
&&\E\Big[\big\|\Delta_sX(t)\big\|_2\big\|\Delta_tX(u)\big\|_2\Big]^2\leq C\big[\E\big[\|\Delta_tX(u)\|_2^4\big]^\frac{1}{2}\big[\E\big[\|\Delta_sX(t)\|_2^4\big]^\frac{1}{2}\nonumber
\\&&\qquad\qquad\leq C \Big(\frac{\lfloor nt_1 \rfloor-\lfloor ns_1\rfloor}{n}\Big)^{H}\Big(\frac{\lfloor nt_2 \rfloor-\lfloor ns_2 \rfloor}{n}\Big)^{H}\nonumber
\\
&&\qquad\qquad\qquad\times\Big(\frac{\lfloor nu_1\rfloor-\lfloor nt_1\rfloor}{n}\Big)^{H}\Big(\frac{\lfloor nu_2\rfloor-\lfloor nt_2\rfloor}{n}\Big)^{H}.
\eeqlb
Hence, if $u_2-s_2\geq \frac{1}{n}$, then
\beqlb\label{da-31}
\Big|\frac{\lfloor nu_2\rfloor-\lfloor ns_2\rfloor}{n}\Big|^{2H}\leq C |(u_2-s_2)|^{2H}.
\eeqlb
Conversely, if $u_2-s_2<\frac{1}{n}$, then either $u_2$ and $t_2$ or $t_2$ and $s_2$ belong to a same subinterval $[\frac{m}{n}, \frac{m+1}{n})$ for some integer $m$.  Hence \eqref{da-31} still holds.
The other term follows a similar discussion.  The proof is now completed.
\qed

Since $X_{n}(t,s)$, $n\in\N$, are
null on the axes,
by using the criterion given by Bickel and Wichura \cite{Bickel}, and Lemma 3.1, we can get the following lemma.
\begin{lem}\label{lem-3}
The sequence $\{X_n(t,s);(t,s)^T\in[0,\;1]^2\}$ is tight in $\mathcal{D}([0,\;1]^2)$.
\end{lem}

Now, in order to prove Theorem  \ref{thm-2}, it suffices to show the following lemma which states that the law of
all possible weak limits is the law of the operator fractional Brownian sheet of Riemann-Liouville type $X$.
\begin{lem}\label{thm-1}
The family of random fields $X_{n}(t, s)$ defined by \eqref{da-30} converges, as $n$ tends to infinity, to the operator
fractional Brownian sheet of Riemann-Liouville type $X$ in the sense of finite-dimensional distributions.
\end{lem}

In order to prove Lemma \ref{thm-1},  we need a technical result.  Before we present this result, we first introduce the following notation.
\beqnn
(t-u)_+^{\frac{D}{2}-\frac{I}{2}}=\big(\tilde{K}_{i,j}(t,u)\big)_{d\times d}
\eeqnn
and
\beqnn
\Big(\frac{\lfloor nt\rfloor}{n}-u\Big)_+^{\frac{D}{2}-\frac{I}{2}}=\big(\tilde{K}^n_{i,j}(t,u)\big)_{d\times d}.
\eeqnn

\begin{lem}\label{da-lem}
For any $(t_{k}, s_{k})^T, (t_{l}, s_{l})^T\in[0, 1]^2$ and $q,m\in\{1,\cdots,d\}$,  we have that
\beqlb\label{da-19}
&n^{4}\sum_{i=1}^{n}\sum_{j=1}^{n}\int_{\frac{i-1}{n}}^{\frac{i}{n}}\int_{\frac{j-1}{n}}^{\frac{j}{n}}
\tilde{K}_{q,m}^n(t_k,u)\tilde{K}_{m,q}^n(s_k,v)dudv\nonumber
\\&\qquad\quad\quad\int_{\frac{i-1}{n}}^{\frac{i}{n}}\int_{\frac{j-1}{n}}^{\frac{j}{n}}
\tilde{K}_{q,m}^n(t_l,u)\tilde{K}_{m,q}^n(s_l,v)dudv(\xi_{i, j,q}^{(n)})^{2}
\eeqlb
converges to
\beqlb\label{da-20}
\int_{0}^{1}\int_0^1 \tilde{K}_{q,m}(t_k,u)\tilde{K}_{m,q}(s_k,v) \tilde{K}_{q,m}(t_l,u)\tilde{K}_{m,q}(s_l,v)dudv,\,\textrm{a.s.}
\eeqlb
as $n$ tends to infinity.
\end{lem}
{\it Proof:}
It is obvious that \eqref{da-19} is equivalent to
\beqlb\label{da-21}
&n^{2}\sum_{i=1}^{n}n\int_{\frac{i-1}{n}}^{\frac{i}{n}}\tilde{K}^n_{q,m}(t_k,u)du\int_{\frac{i-1}{n}}^{\frac{i}{n}}\tilde{K}_{q,m}^n(t_l,u)du\nonumber\\
&\qquad\qquad\cdot\sum_{j=1}^{n}n\int_{\frac{j-1}{n}}^{\frac{j}{n}}\tilde{K}^n_{m,q}(s_k,v)dv\int_{\frac{j-1}{n}}^{\frac{j}{n}}\tilde{K}^n_{m,q}(s_l,v)dv(\xi_{i, j,q}^{(n)})^{2}.
\eeqlb
By using the same method as the proof of Lemma 8 in Dai,  Hu and Lee\cite{Hong Shuai Dai 1}, we can prove the lemma.
\qed

Next, we prove Lemma \ref{thm-1}.

\noindent{\bf Proof of Lemma \ref{thm-1}.} Let $a_{1}, ..., a_{Q}\in\mathbb{R}$ and $(t_{1}, s_{1})^T, ..., (t_{Q}, s_{Q})^T\in[0,\; 1]^2$. Next,
we prove that the random vector
\begin{align*}
Y_{n}=\sum_{k=1}^{Q}a_{k}X_{n}(t_{k}, s_{k})
\end{align*}
converges in distribution, as $n$ tends to infinity, to the Gaussian random vector
\begin{align*}
\tilde{X}=\sum_{k=1}^{Q}a_{k}X(t_{k}, s_{k}).
\end{align*}
By the well-known Cram\'er-Wold device,  see Whitt \cite{Whitt} for example, in order to prove the above statement, we only need to show that as $n\to\infty$
\beqlb\label{da-7}
bY_n\stackrel{\mathscr{D}}{\to} b\tilde{X},
\eeqlb
where $b=(b_1,b_2,\cdots,b_d)$ and $\stackrel{\mathscr{D}}{\to}$ denotes convergence in distribution.

For conciseness of the paper, let
\beqnn
(t-u)_+^{\frac{D}{2}-\frac{I}{2}}(s-v)_{+}^{\frac{D}{2}-\frac{I}{2}}=K(t,s,u,v)=\big(K_1(t,s,u,v),\cdots, K_d(t,s,u,v)\big)^T,
\eeqnn
where $$K_j(t,s,u,v)=\big(K_{j,1}(t,s,u,v),\cdots, K_{j,d}(t,s,u,v)\big).$$   Then, we have
\beqnn
bY_n&&=\sum_{q=1}^{d}\sum_{k=1}^{Q}\sum_{i=1}^{\lfloor nt \rfloor}\sum_{j=1}^{\lfloor ns \rfloor} n^2 \int_{\frac{i-1}{n}}^{\frac{i}{n}}\int_{\frac{j-1}{n}}^{\frac{j}{n}} a_kb_qK_q (\frac{\lfloor nt_k\rfloor}{n},\frac{\lfloor ns_k\rfloor}{n},u.v)\eta_{i,j}^{(n)}dudv
\\&&=\sum_{m=1}^d\sum_{q=1}^d\sum_{k=1}^Q\sum_{i=1}^{\lfloor nt\rfloor}\sum_{j=1}^{\lfloor ns\rfloor} n^2 \int_{\frac{i-1}{n}}^{\frac{i}{n}}\int_{\frac{j-1}{n}}^{\frac{j}{n}}a_kb_qK_{q,m}(\frac{\lfloor nt_k\rfloor}{n},\frac{\lfloor ns_k\rfloor}{n},u,v)\xi_{i,j,m}^{(n)}dudv,
\eeqnn
and
\beqnn
b\tilde{X}=\sum_{m=1}^d\sum_{q=1}^d\sum_{k=1}^Q \int_{0}^{1}\int_{0}^{1}a_kb_qK_{q,m}(t_k,s_k,u,v)B^{m}(du,dv).
\eeqnn
Since $\xi_{i,j,m}^{(n)}$, $m=1,\cdots,d,$ are independent, in order to prove \eqref{da-7}, we only need to show
\beqlb\label{da-9}
&&\sum_{q=1}^d\sum_{k=1}^Q\sum_{i=1}^{\lfloor nt\rfloor}\sum_{j=1}^{\lfloor ns\rfloor} n^2 \int_{\frac{i-1}{n}}^{\frac{i}{n}}\int_{\frac{j-1}{n}}^{\frac{j}{n}}a_kb_qK_{q,m}(\frac{\lfloor nt_k\rfloor}{n},\frac{\lfloor ns_k\rfloor}{n},u,v)\xi_{i,j,m}^{(n)}dudv\stackrel{\mathscr{D}}{\to}\nonumber
\\&&\qquad\qquad\qquad\qquad\sum_{q=1}^d\sum_{k=1}^Q \int_{0}^{1}\int_{0}^{1}a_kb_qK_{q,m}(t_k,s_k,u,v)B^{m}(du,dv).
\eeqlb

For convenience, we introduce the following notation.
\beqnn
Y^{(n)}_{i,j}=\sum_{q=1}^d\sum_{k=1}^Q n^2\int_{\frac{i-1}{n}}^{\frac{i}{n}}\int_{\frac{j-1}{n}}^{\frac{j}{n}}a_kb_qK_{q,m}(\frac{\lfloor nt_k\rfloor}{n},\frac{\lfloor ns_k\rfloor}{n},u,v)\xi_{i,j,m}^{(n)}dudv.
\eeqnn
Then, \eqref{da-9} can be rewritten as
\beqlb\label{da-10}
\sum_{i=1}^{\lfloor nt\rfloor}\sum_{j=1}^{\lfloor ns\rfloor}Y^{(n)}_{i,j}\stackrel{\mathscr{D}}{\to} \sum_{q=1}^d\sum_{k=1}^Q \int_{0}^{1}\int_{0}^{1}a_kb_qK_{q,m}(t_k,s_k,u,v)B^{m}(du,dv).
\eeqlb

Inspired by Wang, Yan and Yu \cite{Zhi Wang 1}, in order to prove \eqref{da-10}, we first prove the following Lindeberg condition
\beqlb\label{da-11}
\lim_{n\rightarrow\infty}\sum_{i=1}^{n}\sum_{j=1}^{n}\E\Big[(Y_{i, j}^{(n)})^{2}1_{\{|Y_{i, j}^{(n)}|>\varepsilon\}}\Big|\mathcal{F}_{i-1, j-1}^{(n)}\Big]=0
\eeqlb
for all $\varepsilon>0$.

In fact, we have that
\beqlb\label{da-12}
&&\Big(n^2\int_{\frac{i-1}{n}}^{\frac{i}{n}}\int_{\frac{j-1}{n}}^{\frac{j}{n}}K_{q,m}(\frac{\lfloor nt_{k}\rfloor}{n},\frac{\lfloor st_{k}\rfloor}{n},u,v )\xi^{(n)}_{i,j,m}dudv\Big)^2\nonumber
 \\&&\qquad\qquad\qquad \leq n^4 \big(\xi_{i,j,m}^{(n)}\big)^2\bigg(\int_{\frac{i-1}{n}}^{\frac{i}{n}}\int_{\frac{j-1}{n}}^{\frac{j}{n}}\big|K_{q,m}(\frac{\lfloor nt_{k}\rfloor}{n},\frac{\lfloor st_{k}\rfloor}{n},u,v )\big|dudv\bigg)^2\nonumber
 \\&&\qquad\qquad\qquad\leq C n^2 (\xi_{i,j,m}^{(n)}\big)^2 \int_{\frac{i-1}{n}}^{\frac{i}{n}}\int_{\frac{j-1}{n}}^{\frac{j}{n}}\big|K_{q,m}(\frac{\lfloor nt_{k}\rfloor}{n},\frac{\lfloor st_{k}\rfloor}{n},u,v )\big|^2dudv.
\eeqlb

It is easy to verify that there exists  some $\delta>0$ with $\lambda_D-\delta>0$ such that
 \beqlb\label{da-13}
\int_{\frac{i-1}{n}}^{\frac{i}{n}} \big\|(t-u)_+^{\frac{D}{2}-\frac{I}{2}}\big\|^2du&&\leq C\int_{\frac{n-1}{n}}^{1}(1-u)_+^{\lambda_D-1-\delta}du,
\eeqlb
since $0<\lambda_D-\delta<1$ and $t\in[0,\,1]$.

Noting the form of $K$,  we get from \eqref{da-12} and \eqref{da-13} that
\beqlb\label{da-14}
&&\Big(n^2\int_{\frac{i-1}{n}}^{\frac{i}{n}}\int_{\frac{j-1}{n}}^{\frac{j}{n}}K_{q,m}(\frac{\lfloor nt_{k}\rfloor}{n},\frac{\lfloor st_{k}\rfloor}{n},u,v )\xi^{(n)}_{i,j,m}dudv\Big)^2\nonumber
 \\&&\qquad\qquad\qquad \leq Cn^2 (\xi_{i,j,m}^{(n)})^2\delta_n^2,
\eeqlb
where
\beqnn
\delta_n=\int_{\frac{n-1}{n}}^{1}(1-u)_+^{\lambda_D-1-\delta}du.
\eeqnn
It follows from \eqref{da-12} and \eqref{da-14} that
\beqlb\label{da-15}
\Big(Y^{(n)}_{i,j}\Big)^2\leq C \sum_{q=1}^{d}\sum_{k=1}^{Q}n^2 a_k^2 b^2_q (\xi_{i,j,m}^{(n)})^2\delta_n^2.
\eeqlb
On the other hand,
\beqlb\label{da-16}
\{|Y_{i, j}^{(n)}|>\varepsilon\}=\{|Y_{i, j}^{(n)}|^{2}>\varepsilon^{2}\}.
\eeqlb
Hence, from  \eqref{da-15} and \eqref{da-16},
\beqlb\label{da-17}
\Big\{|Y_{i, j}^{(n)}|>\varepsilon\Big\}\subseteq\Big\{ C n^2 (\xi_{i,j,m}^{(n)})^2\delta_n^2>\epsilon^2\Big\}.
\eeqlb
Consequently,
\beqlb\label{da-18}
\E\bigg[(Y_{i, j}^{(n)})^{2}1_{\{|Y_{i, j}^{(n)}|>\varepsilon\}}\Big|\mathcal{F}_{i-1, j-1}^{(n)}\bigg]&&\leq C\E\bigg[n^{2}(\xi_{i, j,m}^{(n)})^{2}\delta_{n}^2
1_{\{Cn^{2}(\xi_{i,j,m}^{(n)})^{2}\delta^{2}_n>\varepsilon^{2}\}}\Big|\mathcal{F}_{i-1, j-1}^{(n)}\bigg]\nonumber\\
&&\leq C\delta_{n}^2\E\bigg[1_{\{Cn^{2}(\xi_{i,j,m}^{(n)})^{2}\delta_{n}^2>\varepsilon^{2}\}}\Big|\mathcal{F}_{i-1, j-1}^{(n)}\bigg]
\eeqlb
for all $i, j=1, 2, ..., n$. Hence, from \eqref{a-1} and \eqref{da-18},
\begin{align*}
&\sum_{i=1}^{n}\sum_{j=1}^{n}\E\bigg[(Y_{i, j}^{(n)})^{2}1_{\{|Y_{i, j}^{(n)}|>\varepsilon\}}\Big|\mathcal{F}_{i-1, j-1}^{(n)}\bigg]\nonumber\\
&\quad\leq\sum_{i=1}^{n}\sum_{j=1}^{n}C\delta_{n}^2\E\bigg[1_{\{Cn^{2}(\xi_{i, j,m}^{(n)})^{2}\delta^2_{n}>\varepsilon^{2}\}}\Big|\mathcal{F}_{i-1, j-1}^{(n)}\bigg]\nonumber\\
&\quad\leq C\delta^2_{n}\sum_{i=1}^{n}\sum_{j=1}^{n}\E[1_{\{C\delta^2_{n}>\varepsilon^{2}\}}]\rightarrow0~~~~(n\rightarrow\infty),
\end{align*}
because $\delta_{n}\rightarrow0$ implies that
$1_{\{C\delta^2_{n}>\varepsilon^{2}\}}=0$ for large enough $n$.

In order to prove \eqref{da-7}, we  also need to show that
\beqlb\label{da-22}
&&\sum_{i=1}^{n}\sum_{j=1}^n\Big[Y^{(n)}_{i,j}\Big]^2\stackrel{\P}{\to}\E\Big[\sum_{q=1}^d\sum_{k=1}^Q \int_{0}^{1}\int_{0}^{1}a_kb_qK_{q,m}(t_k,s_k,u,v)B^{m}(du,dv)\Big],
\eeqlb
where $\stackrel{\P}{\to}$ denotes convergence in probability.
For convenience,  we define
\beqnn
\tilde{B}^m(t,s,u,v)=\sum_{q=1}^db_qK_{q,m}(t,s,u,v).
\eeqnn
Note that the right-hand side  of \eqref{da-22} is equivalent to
\beqlb\label{da-23}
\sum_{i,j=1}^Qa_ia_j\int_0^1\int_0^1\tilde{B}^m(t_i,s_i,u,v)\tilde{B}^m(t_j,t_j,u,v)dudv.
\eeqlb
Next, we look at the left-hand side of \eqref{da-22}. In fact, we have
\beqlb\label{da-24}
Y^{(n)}_{i,j}=\sum_{k=1}^Qa_k\int_{\frac{i-1}{n}}^{\frac{i}{n}}\int_{\frac{j-1}{n}}^{\frac{j}{n}}\tilde{B}^m(\frac{\lfloor nt_k\rfloor}{n},\frac{\lfloor ns_k\rfloor}{n},u,v)\xi^{(n)}_{i,j,m}dudv.
\eeqlb
Hence,
\beqlb\label{da-25}
\Big(Y^{(n)}_{i,j}\Big)^2=&&\big(\xi^{(n)}_{i,j,m}\big)^2\sum_{k,l=1}^{Q}a_ka_l\int_{\frac{i-1}{n}}^{\frac{i}{n}}\int_{\frac{j-1}{n}}^{\frac{j}{n}}\tilde{B}^m(\frac{\lfloor nt_k\rfloor}{n},\frac{\lfloor ns_k\rfloor}{n},u,v)dudv \nonumber
\\&&\cdot \int_{\frac{i-1}{n}}^{\frac{i}{n}}\int_{\frac{j-1}{n}}^{\frac{j}{n}}\tilde{B}^m(\frac{\lfloor nt_l\rfloor}{n},\frac{\lfloor ns_l\rfloor}{n},u,v)dudv.
\eeqlb
Here, we  point out that the entry $K_{q,m}(t,s,u,v)$ takes the form of
\beqnn
\sum_{i=1}^d \tilde{K}_{q,i}(t,u)\tilde{K}_{q,i}(s,v)\tilde{K}_{i,m}(t,u)\tilde{K}_{i,m}(s,v).
\eeqnn
Hence, it follows from Lemma \ref{da-lem} and \eqref{da-23}-\eqref{da-25} that \eqref{da-22} holds.

From the above arguments, we can easily get that the lemma holds.\qed

Now, we prove Theorem \ref{thm-2}.

\noindent{\bf Proof of Theorem \ref{thm-2}:} Theorem \ref{thm-2} is a direct consequence of Lemmas \ref{lem-3} and \ref{thm-1}, because
tightness and the convergence of finite dimensional distributions imply weak convergence (see Bickel and Wichura \cite{Bickel}).
\section{Final Note}
In this work,  based on the fractional Brownian motion of Riemann-Liouville type, we introduce the operator fractional Brownian  sheet of Riemann-Liouville type $X$ and present an approximation to it via  martingale differences.  In Definition \ref{da-35},  $\lambda_D$ and $\Lambda_D$ are assumed to be at $(0,\,1)$. In fact,  if we only want to define a random field $X$, $\lambda_D$ and $\Lambda_D$ can be at a larger interval than $(0,\,1)$. However, in this paper, we also need the random field $X$ to   enjoy some nice properties.
 It follows from Mason and Xiao \cite{Mason and Xiao} that for an operator self-similar random filed $\{\hat{X}(t,s)\}$ with exponent $\hat{D}$, if $\lambda_{\hat{D}}>0$, then $\hat{X}(0,0)=(0,\cdots,0)^T$ a.s. Furthermore, if $\hat{X}(1,0)$ is proper and $\E[\|\hat{X}(1,0)\|_2]<\infty$,  then $\Lambda_{\hat{D}}\leq 1$.  In this  paper,  the operator fractional Brownian sheet of Riemann-Liouville type $\{X(t,s)\}$ is assumed to be proper for  $(t,s)^T=(1,0)^T$, and $X(0,0)=(0,\cdots,0)^T$ a.s. Hence, we assume   $0<\lambda_D,\Lambda_D<1$ in  \eqref{da-8}.

On the other hand, we get from Ayache, L\`{e}ger and Pontier \cite{ALP2002} that, in the one dimensional case $(d=1)$, a fractional Brownian sheet $\{W^{\alpha,\beta}(t,s)\}$ with two parameters $\alpha, \beta\in(0,\;1)$ can be defined as
\beqlb\label{da-32}
W^{\alpha,\beta}(t,s)=\int_{\R^2}f_{\alpha}(t,u)f_{\beta}(s,v)\tilde{B}(dv, du),
\eeqlb
where $f_{H}(t,u)=(t-u)_+^{H-\frac{1}{2}}-(-u)_+^{H-\frac{1}{2}}$. Hence in the one-dimensional case ($d=1$), $X$ defined by \eqref{da-8} is a special kind of fractional Brownian sheets (with $\alpha=\beta$) of Riemann-Liouville type. In fact, inspired by \eqref{da-32},  one could like to define the operator fractional Brownian sheet of Riemann-Liouville type $\hat{X}=\{\hat{X}(t,s);(t,s)^T\in\R_+^2\}$  by
\beqnn
\hat{X}(t,s)=\int_{0}^\infty\int_{0}^\infty (t-u)_+^{\frac{D}{2}-\frac{I}{2}}(s-v)_+^{\frac{\hat{D}}{2}-\frac{I}{2}}B(du,dv),
\eeqnn
where $\hat{D}$ is a linear operator on $\R^d$ with $0<\lambda_{\hat{D}},\Lambda_{\hat{D}}<1$ . It is easy to verify that the random field $\hat{X}$ is well defined.  However,  in such case,  we can not get Theorem \ref{thm-2} according to our method.

\noindent{\bf Acknowledgments}\; We thank the referee and the editor for their time and comments.   Dai was supported by the  National Natural Science Foundation of China (No.11361007),  the Shandong Natural Science Foundation (No.  ZR2014AM021) and  the Fostering Project of Dominant Discipline and Talent Team of Shandong Province Higher Education Institutions.  Shen was supported by the National Natural Science Foundation of China (No.11271020), the Distinguished Young Scholars Foundation of Anhui Province (No. 1608085J06), and the Top Talent Project of University Discipline (Speciality) (No. gxbjZD03)


\end{document}